\documentclass[twocolumn,twoside]{IEEEtran}

\ifCLASSINFOpdf

\else

\fi

\ifCLASSOPTIONcompsoc
\usepackage[caption=false, font=normalsize, labelfont=sf, textfont=sf]{subfig}
\else
\usepackage[caption=false, font=normalsize]{subfig}
\fi
\usepackage{lipsum}%
\usepackage[dvipsnames]{xcolor}
\usepackage{algorithm,algorithmic}
\usepackage{balance}
\usepackage{multicol}   
\usepackage{cite}
\usepackage{gensymb}
\usepackage{multirow}
\usepackage{graphics}
\usepackage{epsfig}
\usepackage{graphicx}
\usepackage{epstopdf}
\usepackage{textcomp}
\usepackage{amsmath}
\usepackage{mathtools}
\usepackage{filecontents}
\usepackage{lipsum,color}
\usepackage{amssymb}
\usepackage{float}
\usepackage{booktabs}
\usepackage{colortbl} 
\usepackage{times} 
\usepackage{amsthm}  

\usepackage{amsfonts}

\theoremstyle{break}
\begin{document}
	\raggedbottom
	\makeatletter
	\def\fps@figure{!t} 
	\def\fps@table{!t}   
	\makeatother

\title{Q-Learning for 3D Coverage in VCSEL-based Optical Wireless Systems}
	
	\author{Hossein~Safi, ~{\it Member,~IEEE}, Rizwana Ahmad,\\~Iman~Tavakkolnia,~{\it Senior Member,~IEEE}, ~and Harald~Haas,~{\it Fellow,~IEEE}
		\thanks{The authors are with the LiFi R\&D centre, Department of Engineering, University of Cambridge, Cambridge, UK, e-mails: \{hs905, ra714, it360, huh21\}@cam.ac.uk.}}

	\maketitle

	\begin{abstract}
		Beam divergence control is a key factor in maintaining reliable coverage in indoor optical wireless communication (OWC) systems as receiver height varies. Conventional systems employ fixed divergence angles, which result in significant coverage degradation due to the non-convex trade-off between optical power concentration and spatial spread. In this paper, we introduce a reinforcement learning (RL)–based framework for dynamic divergence adaptation in vertical-cavity surface-emitting laser (VCSEL)–based OWC networks. By continuously interacting with the environment, the RL agent autonomously learns a near-optimal mapping between receiver height and beam divergence, thereby eliminating the need for analytical modeling or computationally intensive exhaustive search. Simulation results demonstrate that the proposed approach achieves up to 92\% coverage at low receiver heights and maintains robust performance under challenging conditions, enabling scalable, real-time, and energy-efficient beam control for dense VCSEL array deployments in next-generation OWC systems.
	\end{abstract}
	\begin{IEEEkeywords}
	Adaptive coverage control, beam divergence optimization, optical wireless communication, reinforcement learning, VCSEL arrays.
	\end{IEEEkeywords}
	\IEEEpeerreviewmaketitle
	
\section{Introduction}

The growing demand for high-capacity, low-latency connectivity in 6G and beyond is driving innovation in wireless technologies. Among the emerging access technologies, optical wireless communication (OWC) has gained considerable attention as a promising complement to radio frequency (RF) systems, offering unlicensed spectrum, enhanced physical-layer security, and the potential to achieve data rates well beyond 100~Gb/s~\cite{6G-NoN1,6G-2}. Recently, vertical-cavity surface-emitting lasers (VCSELs) have been proposed as key enablers for high-performance OWC access points (APs), thanks to their high modulation bandwidth, energy efficiency, and scalable compact array design~\cite{VCSELbook, LiFi2}. VCSEL arrays can generate multiple narrow and highly directional optical beams that simultaneously serve different regions within an indoor environment, enabling aggregate data rates exceeding the terabit-per-second range \cite{LiFi2}. Their circular beam symmetry and wafer-level manufacturability make them ideal for forming structured dense grid arrays with independently addressable emitters \cite{VCSEL1}.

Despite these advantages, practical VCSEL-based OWC systems still face several critical challenges related to spatial dynamics and beam geometry. Existing array designs typically rely on fixed optical configurations, such as lenses with constant focal length, that are optimized for a specific receiver plane height~\cite{VCSEL1,VCSEL2,VCSEL3}. As a result, when users move vertically or when the receiver plane changes, the overlap and divergence of adjacent beams become suboptimal, causing severe coverage degradation, inter-beam interference, and signal power loss. Furthermore, the number of active VCSELs and their beamwidths are often chosen to guarantee coverage at a single height, leading to unnecessary energy expenditure and reduced efficiency at other heights. These issues are exacerbated in dense deployments where the beam divergence must balance two competing objectives: maximizing area coverage while minimizing optical interference.

To address these limitations, dynamic beam divergence control has emerged as a key design objective for future OWC networks. Recent advances in adaptive optical elements, such as electrically tunable liquid crystal lenses and variable-focus optics, now make it feasible to adjust the divergence angle of each VCSEL beam in real time~\cite{vlenses,vlenses2}. However, determining the optimal divergence configuration remains a non-convex problem with strong coupling between transmitter geometry, propagation loss, and spatial interference. Traditional analytical or gradient-based optimization methods are not suitable for such settings, as the objective function is piecewise, non-differentiable, and highly environment-dependent \cite{safi3d}. 

Motivated by these challenges, this paper introduces an artificial intelligence (AI)–driven approach that leverages Q-learning to autonomously optimize beam divergence angles and enable dynamic three-dimensional (3D) coverage in VCSEL-based OWC systems. \textcolor{black}{More precisely, we formulate the beam selection and activation problem using reinforcement learning (RL), a class of sequential decision-making methods in which an agent learns an optimal policy through interaction with the environment to maximize a long-term reward.} Unlike static configurations or pre-computed lookup tables, the proposed Q-learning framework learns an adaptive state–action policy that directly maps environmental conditions (e.g., receiver height) to the optimal divergence setting. This model-free formulation eliminates the need for explicit channel modeling and enables near real-time adaptation with lower inference complexity compared to exhaustive or iterative optimization methods.

The effectiveness and practicality of the proposed approach are validated through comprehensive simulations, demonstrating its ability to achieve near-optimal coverage and efficient convergence across diverse receiver heights and deployment scenarios. More precisely, simulation results demonstrate that the proposed approach achieves near-optimal coverage performance, maintaining up to 92\% coverage at low receiver heights and 50\% under challenging conditions, while reducing the online computational cost by several orders of magnitude compared to exhaustive optimization. These results confirm that AI-driven learning can transform beam management in optical wireless systems, paving the way toward intelligent, scalable, and self-optimizing VCSEL access networks for future 6G environments.

The remainder of this paper is organized as follows. Section~\ref{sys-model} describes the system model and beam propagation characteristics. Section~\ref{sec:optimization} formulates the beam divergence optimization problem and motivates the need for adaptive learning. Section~\ref{sec:rl_optimization} discusses the AI-based learning framework. Section~\ref{sec:results} presents the simulation setup and numerical results. Finally, Section~\ref{Conclusion} concludes the paper and outlines future research directions.
	\section{System Model}
	\label{sys-model}
	
	\begin{figure}
		\begin{center} 
			\includegraphics[width=3.1 in]{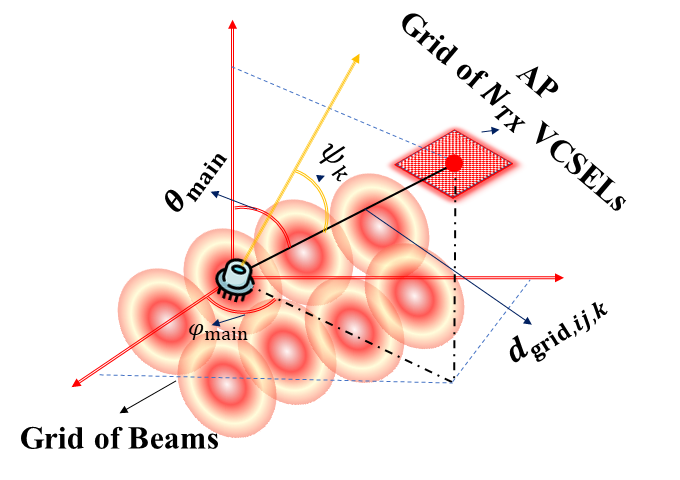}
			\caption{Configuration of VCSEL-Based Downlink Transmission}
			\label{sysmodel}
		\end{center}
	\end{figure}

	
	\subsection{Geometric Configuration and Transmitter Array}
	
	The system operates within a rectangular indoor environment defined by the three-dimensional space $L \times W \times H$. A centralized Access Point (AP) is installed at the center of the ceiling and supports downlink communication. 	The AP employs a structured grid of $N_{\text{TX}}$ VCSELs \cite{VCSEL2,safi3d}, thereby forming a corresponding grid of beams at the receiver plane, as illustrated in Fig.~\ref{sysmodel}. 
	The $k$-th individual VCSEL beam is actively steered from its source position $(x_{\text{source}, k}, y_{\text{source}, k}, z_{\text{source}, k})$. The main beam direction is thus defined by the polar angle $\theta_{\text{main}}$ and the azimuthal angle $\phi_{\text{main}}$, forming the beam direction vector $\mathbf{d}_{\text{beam}, k}$ as
	\begin{equation}
		\mathbf{d}_{\text{beam}, k}=\begin{bmatrix} \sin(\theta_{\text{main}}) \cos(\phi_{\text{main}}) \\ \sin(\theta_{\text{main}}) \sin(\phi_{\text{main}}) \\ \cos(\theta_{\text{main}}) \end{bmatrix}.
	\end{equation}
	
	The area of interest is the receiver plane, which is discretized into grid points $(x_i, y_j)$ with resolutions $\Delta_x$ and $\Delta_y$. The spatial vector from the $k$-th transmitter to the grid point $(i, j)$ is defined as
	\begin{equation}
		\mathbf{d}_{\text{grid}, ij, k}=\begin{bmatrix} x_{i}-x_{\text{source}, k} \\ y_{j}-y_{\text{source}, k} \\ z_{ij}-z_{\text{source}, k} \end{bmatrix},
	\end{equation}
	where $d_{ij, k} = ||\mathbf{d}_{\text{grid}, ij, k}||$ is the Euclidean distance. This geometric formulation establishes the basis for characterizing beam propagation, where the received irradiance at each grid point is determined by the divergence and angular alignment of the Gaussian beam.
	\vspace{-3mm}
	
		\subsection{Gaussian Beam Propagation and Received Irradiance}
	Let $\theta_{\text{diff},ij,k}$ denote the angular deviation between the beam axis $\mathbf{d}_{\text{beam},k}$ and the direction vector toward grid point $\mathbf{d}_{\text{grid},ij,k}$, defined as
	\begin{equation}
		\theta_{\text{diff},ij,k}
		= \arccos\!\left(
		\frac{\mathbf{d}_{\text{beam},k}\!\cdot\!\mathbf{d}_{\text{grid},ij,k}}
		{\|\mathbf{d}_{\text{beam},k}\|\,\|\mathbf{d}_{\text{grid},ij,k}\|}
		\right).
	\end{equation}
	The beam intensity at the receiver plane follows Gaussian beam propagation, with range-dependent beam radius
	$w(z)$. Accordingly, the received irradiance at grid point $(i,j)$ from transmitter $k$ is
	\begin{equation}
		\label{eq:irradiance_wz}
		I_{\text{grid},ij,k}
		= \frac{2P_t}{\pi\,w^2(d_{ij,k})}\,
		\exp\!\left[
		-2\left(\frac{d_{ij,k}\,\theta_{\text{diff},ij,k}}{w(d_{ij,k})}\right)^{\!2}
		\right],
	\end{equation}
	where $P_t$ is the transmitted optical power.  Here, $w(d_{ij,k})$ captures the far-field expansion with distance, while the exponential term models the angular decay of the Gaussian profile. Using $w(d_{ij,k})\!\approx\! d_{ij,k}\,\theta_{\text{divergence}}$ \cite{VCSEL1}, the irradiance in \eqref{eq:irradiance_wz} simplifies to
	\begin{equation}
		\label{eq:irradiance_wz2}
		I_{\text{grid},ij,k}
		\!\approx\! \frac{2P_t}{\pi[\theta_{\text{divergence}}d_{ij,k}]^{2}}
		\exp\!\left[\!-\!2\left(\frac{\theta_{\text{diff},ij,k}}{\theta_{\text{divergence}}}\right)^{\!2}\right]\!.
	\end{equation}
	As shown in (\ref{eq:irradiance_wz2}), the divergence angle governs spatial coverage. Increasing it broadens the beam and coverage area but also raises interference due to power leakage. Thus, optimization is essential to balance coverage and interference.

	\vspace{-2mm}
	\subsection{Receiver Front-End Model and Final Received Power}
	 To obtain the final collected optical power, the receiver front-end model must be applied. The receiver is assumed to be mounted horizontally, giving its normal vector an angle of incidence $\Psi_{k}$ relative to the incoming light ray. The total received power from the $k$-th transmitter is calculated as \cite{VCSEL1}
	\begin{equation}
		P_{\text{rec}, k} = I_{\text{grid}, k} \cdot A_d \cdot \cos(\Psi_{k}) \cdot G_{\text{conc}} \cdot \Pi\!\left(\frac{\Psi_{k}}{\Psi_{\text{FoV}}}\right),
	\end{equation}
	where $A_d$ is the effective photodetector area, $G_{\text{conc}}$ is the optical concentrator gain, and $\Psi_{\text{FoV}}$ is the receiver field-of-view (FoV). The factor $\cos(\Psi_{k})$ captures geometric loss, while the gate function $\Pi(\cdot)$ ensures that only rays within the receiver FoV are collected.
	
	\section{Optimization Framework for Divergence Angle Control}
	\label{sec:optimization}
	
	\subsection{Problem Formulation}
	
	As discussed, the objective is to determine the optimal divergence angle, denoted as $\theta_{\text{opt}}$, that maximizes spatial coverage under a given environmental state $s$. The state $s$ is primarily determined by the receiver height $h_r$, although it may also incorporate additional conditions of the link setup.
	The coverage percentage, denoted as $C$, is defined as the fraction of grid points that satisfy the minimum signal-to-interference-plus-noise ratio (SINR) threshold $\Gamma_{\text{th}}$. The optimization problem can therefore be expressed as
	\begin{equation}
		\label{opt-problem}
		\underset{\theta_{\text{divergence}} \in \mathcal{A}}{\text{maximize}} \quad C(\theta_{\text{divergence}} | s),
	\end{equation}
	where $\mathcal{A}$ is the feasible set of divergence angles constrained by the optical hardware\footnote{In practice, the divergence angle can be dynamically tuned through optical subsystems integrated at the transmitter side, such as liquid crystal lenses, variable-focus optics, or adaptive beam-shaping elements, which enable fine-grained control of the emitted beam profile in real time\cite{vlenses,vlenses2,liquidlens}.}. The coverage function is given by
	\begin{equation}
		C(\theta_{\text{divergence}} | s) = \frac{1}{N_{\text{grid}}} \sum_{i,j} \mathbb{I}\!\left( \text{SINR}_{i,j}(\theta_{\text{divergence}}) \geq \Gamma_{\text{th}} \right),
	\end{equation}
	with $N_{\text{grid}}$ denoting the number of discrete grid locations and $\mathbb{I}(\cdot)$ the indicator function.
	
	At each grid point $(i, j)$, the SINR is defined as
	\begin{equation}
		\label{eq:SINR}
		\begin{split}
			\text{SINR}_{ij}(\theta_{\text{divergence}}) &=
			\frac{P_{\text{rec},ij,k^*}(\theta_{\text{divergence}})}
			{\sum_{l \neq k^*} P_{\text{rec},ij,l}(\theta_{\text{divergence}}) + N_0}, \\
			k^* &= \arg\max_k P_{\text{rec},ij,k}(\theta_{\text{divergence}}).
		\end{split}
	\end{equation}
	where $P_{\text{rec}, k, i, j}$ is the received power from transmitter $k$ at location $(i, j)$ and $k^{*}$ denotes the index of the dominant beam. The denominator captures both residual interference and constant noise power $N_0$. The optimization problem is thus to identify $\theta_{\text{opt}}$ such that $C$ is maximized over the service region.

	\subsection{Analytical Challenges and Dynamic Nature}
	
	The optimization problem in (\ref{opt-problem}) is intractable using conventional methods due to structural and environmental complexity. The coverage function $C(\theta_{\text{divergence}})$ is non-linear and non-convex, and the maximization operator $\max_k\{P_{\text{rec},k}\}$ introduces discontinuities by abruptly switching the serving beam index. Additionally, the indicator function $\mathbb{I}(\cdot)$ creates a piecewise-constant, non-differentiable objective. These properties violate the assumptions required for gradient-based optimization methods and such methods cannot be applied in this setting \cite{boyd2004convex}.
	More precisely, let $f(\theta_{\text{divergence}}) = \text{SINR}_{i,j}(\theta_{\text{divergence}}) - \Gamma_{\text{th}}$. The coverage function can then be written as 
	\begin{equation}
		C(\theta_{\text{divergence}}) = \frac{1}{N_{\text{grid}}} \sum_{i,j} \mathbb{I}(f(\theta_{\text{divergence}}) \geq 0).
	\end{equation}
	Since $\mathbb{I}(\cdot)$ is discontinuous at zero and flat elsewhere, the gradient $\partial C/\partial \theta_{\text{divergence}}$ does not exist in the classical sense, and subgradient approximations yield no useful direction for optimization. This creates a situation similar to mixed-integer programming problems, which are known to be NP-hard \cite{boyd2004convex}. 
	
	Beyond the mathematical intractability, the problem is inherently dynamic. The optimal divergence angle $\theta_{\text{opt}}$ depends on the instantaneous environmental state $s$. As the receiver height changes, the geometry of the link evolves and a previously optimal setting may become sub-optimal. A static divergence angle optimized for a single operating condition leads to significant coverage degradation at other conditions. Similarly, reliance on pre-computed look-up tables is insufficient in environments with fast temporal variation. 	Consequently, the optimization problem exhibits both non-convex structural properties and non-stationary environmental dynamics, motivating the need for adaptive, learning-based control methods rather than closed-form or static solutions.
	\section{Reinforcement Learning Formulation}
	\label{sec:rl_optimization}
	
To address the dual challenges of analytical intractability and dynamic environmental variability discussed in the previous section, we adopt an RL framework. This approach reformulates the divergence control task as a sequential decision-making process, enabling the transmitter to learn an adaptive policy that maps the instantaneous state (receiver height) to the optimal divergence angle. Through iterative interaction with the custom-built environment, the RL agent continuously refines its control strategy, converging to a policy that maximizes the expected cumulative reward derived from coverage performance. This RL-based approach provides real-time adaptability by learning through direct interaction with the environment, allowing the system to optimize performance continuously without relying on explicit analytical models or prior channel knowledge. \textcolor{black}{The beam selection problem is formulated as a finite-state Markov decision process with a discrete action space corresponding to the available VCSEL beams. Under these conditions, tabular Q-learning provides a natural and efficient solution framework. It offers low computational complexity, stable convergence properties, and direct interpretability of the learned action–value function.}

\subsection{Markov Decision Process Formulation}

The divergence optimization task is modelled as a finite Markov Decision Process (MDP) defined by the tuple
\begin{equation}
	\mathcal{M} = (\mathcal{S}, \mathcal{A}, P, R, \gamma),
\end{equation}
where each element represents a key component of the decision-making environment as follows:

\begin{itemize}
		\item \textbf{State space:} $\mathcal{S} = \{s_1, s_2, \dots, s_{N_S}\}$ represents the set of observable environmental states. In this work, each state corresponds to a discrete receiver height. This discretization captures representative operating levels within a typical indoor environment, where users may hold or place devices at table, chest, or standing height. The assumption of discrete states simplifies the learning process while preserving the key spatial variability of the link geometry. The framework, however, can be readily extended to incorporate continuous height values or additional contextual parameters.

	\item \textbf{Action space:} $\mathcal{A}=\{\theta_1,\dots,\theta_{N_A}\}$ denotes the feasible divergence angles. The range and resolution of $\mathcal{A}$ are determined by the optical design and actuation limits of the VCSELs in the array, which bound the extent and granularity of real-time divergence adjustment.

	\item \textbf{Transition probability:} $P(s'|s,a)$ defines the probability of transitioning to state $s'$ after action $a$ in state $s$. 
	We assume that the VCSEL divergence control action $a$ does not influence this transition, as user mobility dynamics are external to the optical control policy \cite{Q-learning}. This decoupling simplifies the MDP’s temporal structure and enables the agent to focus on maximizing the instantaneous coverage reward $R(s,a)$ at the observed state $s$, thereby aligning with the objective of determining the optimal static divergence angle per discrete height.

	\item \textbf{Reward function:} $R(s, a)$ quantifies the system performance after applying action $a$ in state $s$. It is defined as
	\begin{equation}
		R(s, a) = 100 - \text{Hole}(s, a),
	\end{equation}
	where $\text{Hole}(s, a)$ represents the percentage of uncovered grid points (coverage holes) resulting from the divergence configuration. A grid point is considered uncovered if its corresponding $\text{SINR}$ falls below the predefined threshold $\Gamma_{\text{th}}$. Thus, maximizing $R(s, a)$  directly maximizes spatial coverage.

	\item \textbf{Discount factor:} $\gamma \in [0,1]$ controls the weighting of future rewards. Given the quasi-static and deterministic state evolution, optimization is dominated by the instantaneous coverage reward, so a moderate discount is adopted to ensure stable learning while preserving sensitivity to per-state performance.

\end{itemize}

Subsequently, the goal of the RL agent at the transmitter side is to learn an optimal policy
\begin{equation}
	\pi^{*}: \mathcal{S} \rightarrow \mathcal{A},
\end{equation}
that maximizes the expected cumulative discounted reward
\begin{equation}
	J(\pi) = \mathbb{E}\!\left[\sum_{t=0}^{\infty} \gamma^{t} R(s_t, a_t)\right].
\end{equation}
This policy enables the transmitter to autonomously select the divergence angle $\theta_{\text{opt}}$ that achieves the highest expected coverage under each observed receiver condition. 

\textcolor{black}{Note that, although the state transitions are not influenced by the agent’s actions, the problem is formulated within an RL framework to capture the sequential decision-making nature of beam configuration under time-varying conditions. In the special case where state transitions are fully independent of actions, the formulation reduces to a contextual bandit problem. However, adopting an RL formulation enables direct extension to more realistic scenarios involving correlated mobility, constrained movement, or temporal coupling between beam selection decisions, without altering the core framework.}

\subsection{Q-Learning Framework}

We employ a model-free Q-learning algorithm to estimate the optimal action-value function $Q^{*}(s, a)$, which represents the expected long-term return of executing action $a$ in state $s$ and following the optimal policy thereafter. The Q-learning update rule is expressed as \cite{Q-learning}
\begin{equation}
	Q(s,a) \leftarrow Q(s,a) + \alpha \big[ R(s,a) + \gamma \max_{a'} Q(s',a') - Q(s,a) \big],
	\label{eq:q_update}
\end{equation}
where $\alpha$ denotes the learning rate. Over repeated episodes, $Q(s,a)$ converges to $Q^{*}(s,a)$, and the optimal policy can be derived as
\begin{equation}
	\pi^{*}(s) = \arg\max_{a} Q^{*}(s,a).
\end{equation}

An $\epsilon$-greedy exploration strategy is also adopted to balance exploration and exploitation. The exploration rate decays exponentially with training episodes
\begin{equation}
	\epsilon = \epsilon_{\min} + (\epsilon_{\max} - \epsilon_{\min}) e^{-\lambda n_{\text{ep}}},
\end{equation}
where $n_{\text{ep}}$ is the episode index and $\lambda$ the decay rate. This ensures that the agent explores broadly during early training and focuses on fine-tuning its policy in later stages. The full procedure is summarized in Algorithm~\ref{alg:q_learning}.

\subsection{Implementation and Convergence Behaviour}

The proposed learning framework is implemented using a custom environment that accurately replicates the geometric configuration, beam propagation characteristics, and receiver front-end model described in Section~\ref{sys-model}. During training, each episode sequentially traverses the discrete receiver states, while the Q-table is iteratively updated based on the coverage reward obtained for each divergence action.
 As detailed in Algorithm~\ref{alg:q_learning}, the agent interacts repeatedly with the custom  environment, applying an $\epsilon$-greedy policy to balance exploration and exploitation. At each step, the observed coverage reward updates the Q-values according to the  equation in~(\ref{eq:q_update}). Through continuous interaction, the Q-table converges to stable estimates of the optimal action–value function, yielding the policy $\pi^{*}$ that maps each environmental state to its corresponding optimal divergence angle
\begin{equation}
	\theta_{\text{opt}} = \pi^{*}(s_{\text{current}}).
\end{equation}

Convergence typically occurs within $10^3$–$10^4$ episodes, depending on exploration settings. Once trained, divergence adaptation is computationally efficient, requiring only a single Q-table lookup. The RL agent remains robust under non-stationary conditions and continuously refines its policy through interaction.

\begin{algorithm}
	\caption{Q-Learning for Adaptive Divergence Control}
	\label{alg:q_learning}
	\begin{algorithmic}[1]
		\REQUIRE State set $\mathcal{S}$, action set $\mathcal{A}$; learning rate $\alpha$; discount $\gamma$; episodes $N_{\text{ep}}$; exploration parameters $(\epsilon_{\max}, \epsilon_{\min}, \lambda)$.
		\STATE Initialize Q-table $Q(s,a)\leftarrow 0,\;\forall s\in\mathcal{S},a\in\mathcal{A}$.
		\STATE Initialize environment $\mathcal{E}$ (\texttt{VCSELEnv}).
		\FOR{$n=1$ to $N_{\text{ep}}$}
		\STATE $s \leftarrow \mathcal{E}.\text{reset}()$ \hfill $\triangleright$ start episode at initial state
		\STATE $\epsilon \leftarrow \epsilon_{\min} + (\epsilon_{\max}-\epsilon_{\min}) e^{-\lambda n}$ \hfill $\triangleright$ episode-wise decay
		\STATE $\text{done} \leftarrow \textbf{false}$
		\WHILE{not $\text{done}$}
		\STATE \hfill $\triangleright$ $\epsilon$-greedy action selection
		\IF{$\text{Uniform}(0,1) < \epsilon$}
		\STATE $a \leftarrow \text{Random}(\mathcal{A})$ \hfill $\triangleright$ explore
		\ELSE
		\STATE $a \leftarrow \arg\max_{a' \in \mathcal{A}} Q(s,a')$ \hfill $\triangleright$ exploit; break ties arbitrarily
		\ENDIF
		\STATE $(s', R, \text{done}) \leftarrow \mathcal{E}.\text{step}(a)$ \hfill $\triangleright$ apply divergence action, observe reward
		\STATE \hfill $\triangleright$ Bellman update (no bootstrapping at terminal)
		\IF{$\text{done}$}
		\STATE $Q(s,a) \leftarrow Q(s,a) + \alpha \big[R - Q(s,a)\big]$
		\ELSE
		\STATE $Q(s,a) \leftarrow Q(s,a) + \alpha \big[R + \gamma \max_{a' \in \mathcal{A}} Q(s',a') - Q(s,a)\big]$
		\ENDIF
		\STATE $s \leftarrow s'$
		\ENDWHILE
		\ENDFOR
		\STATE \textbf{Output:} $\pi^{*}(s) = \arg\max_{a \in \mathcal{A}} Q(s,a)$.
	\end{algorithmic}
\end{algorithm}

\section{Simulation Results and Discussion}
\label{sec:results}
\begin{figure*}
	\centering
	\subfloat[] {\includegraphics[width=1.72 in]{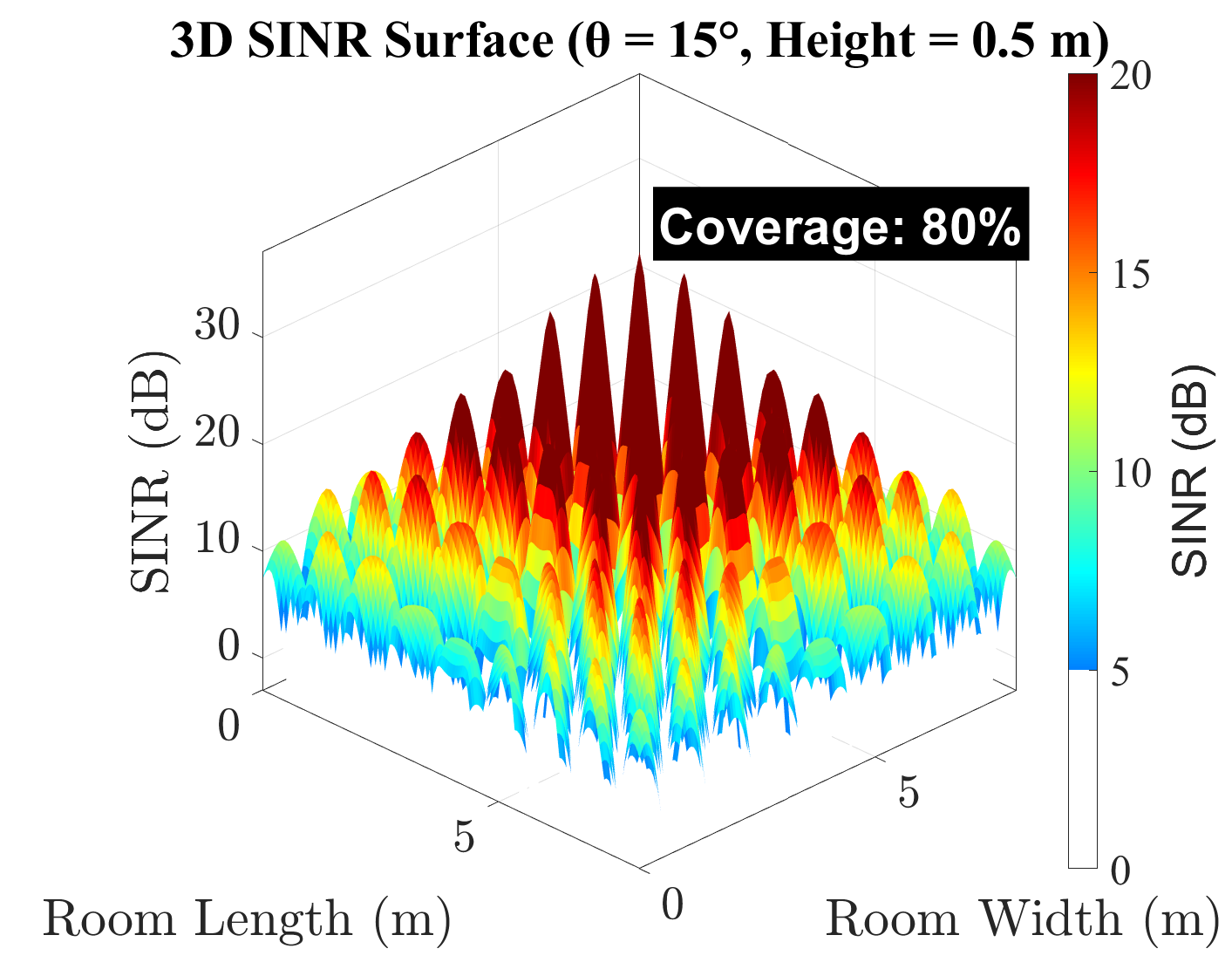}
		\label{cf1}
	}
	\hfill
	\subfloat[] {\includegraphics[width=1.72 in]{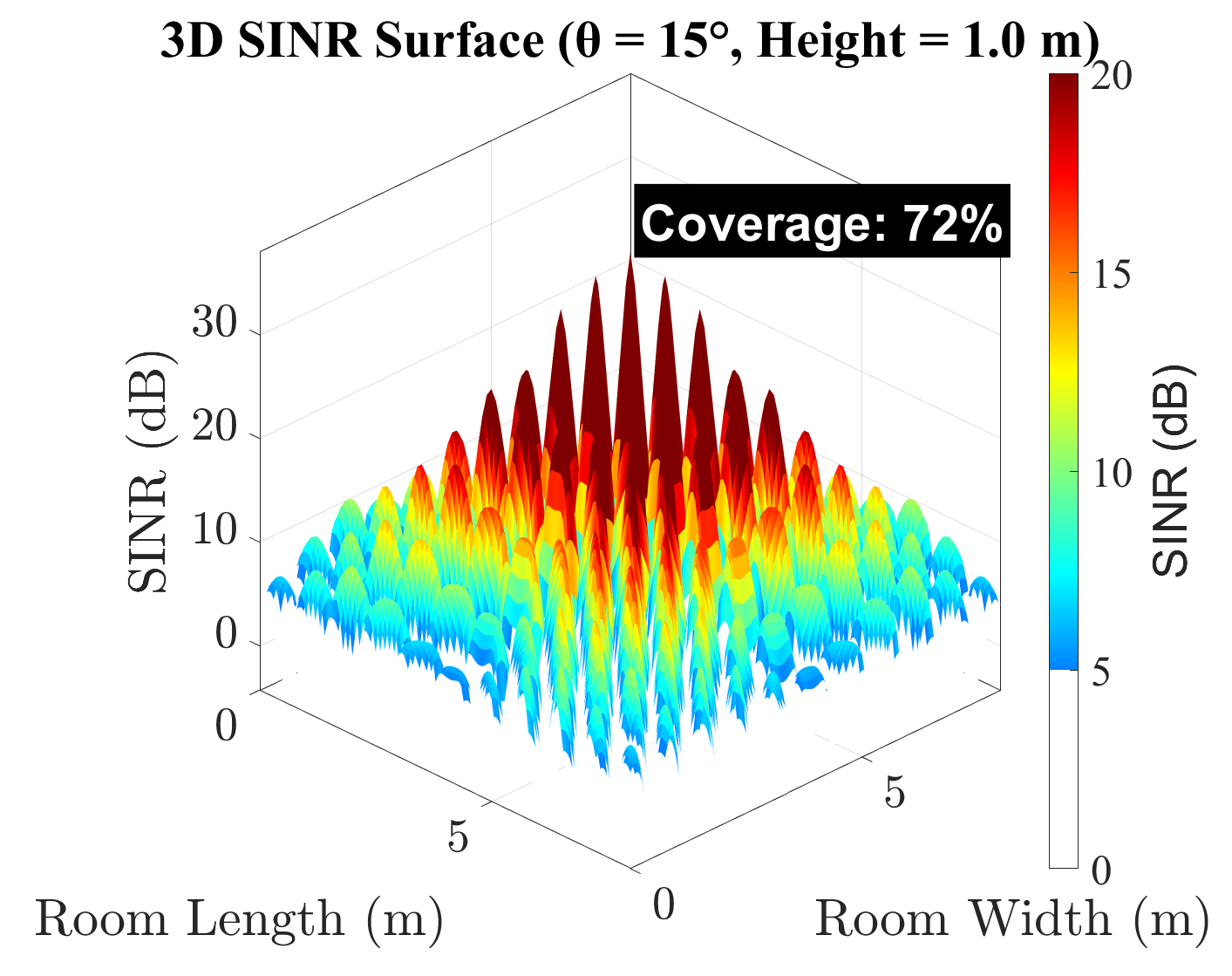}
		\label{cf2}
	}
	\hfill
	\subfloat[] {\includegraphics[width=1.72 in]{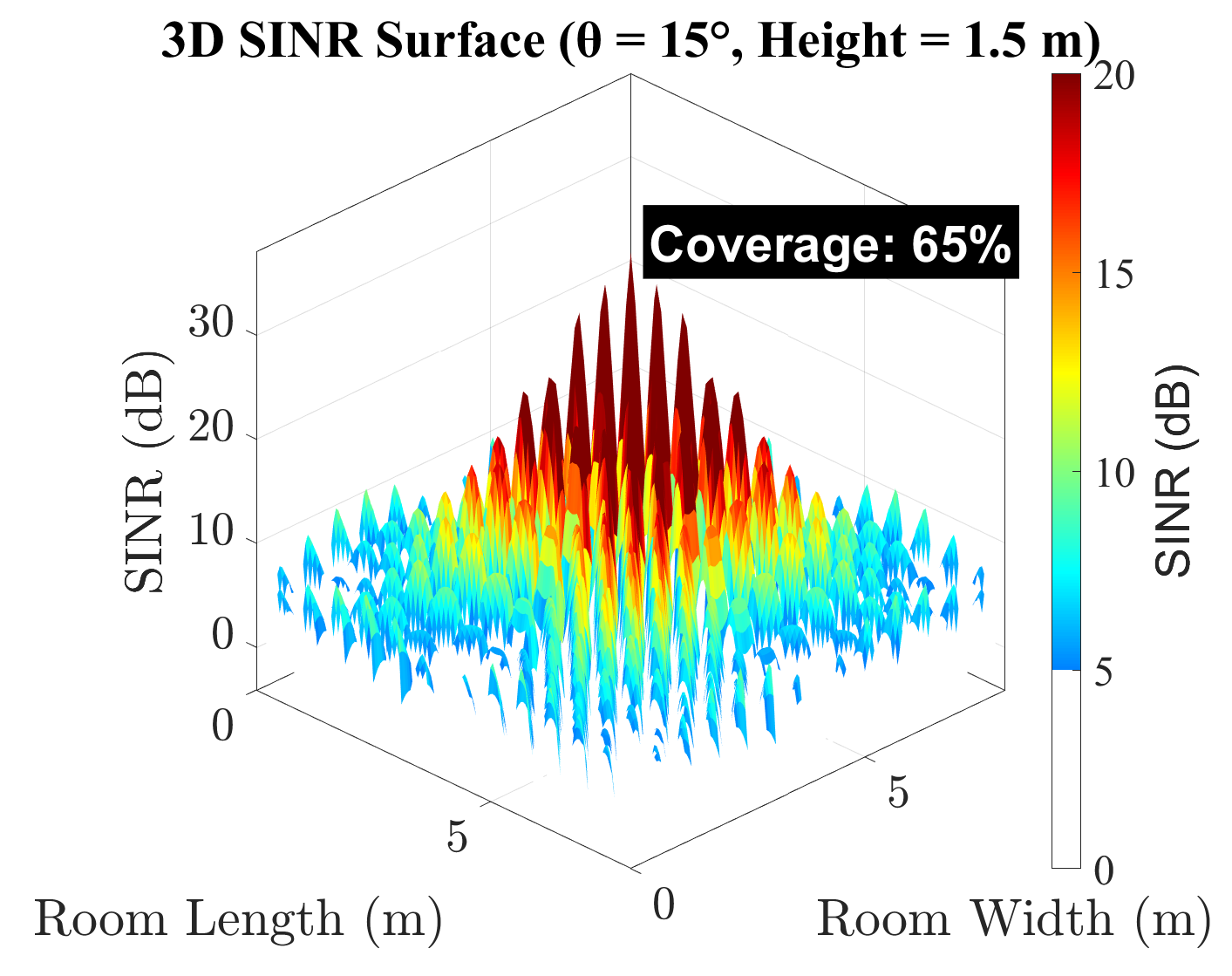}
		\label{cf3}
	}
	\hfill
	\subfloat[] {\includegraphics[width=1.72 in]{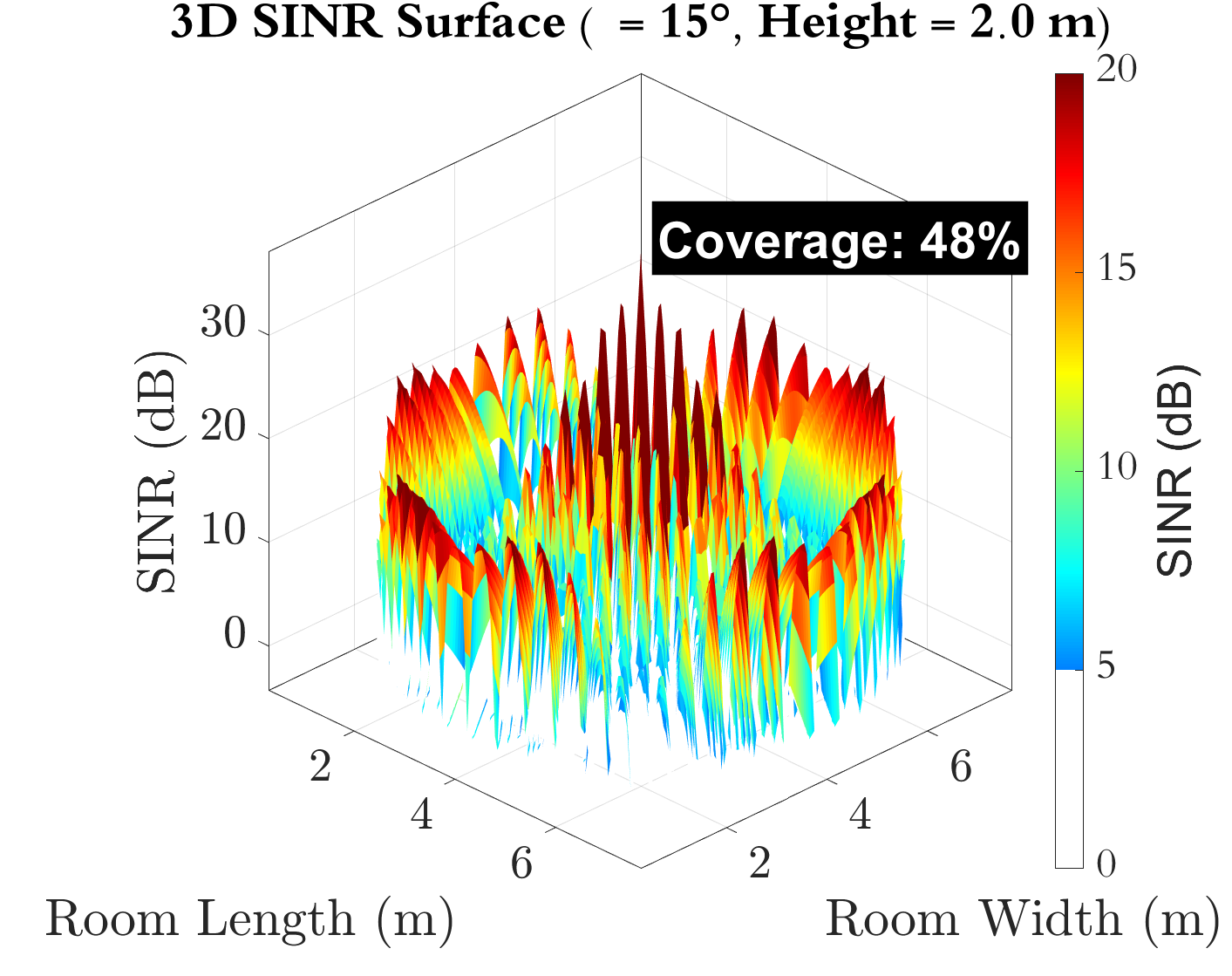}
		\label{cf4}
	}
	\caption{Coverage performance for $\theta_{\text{divergence}} = 15^\circ$ at different heights. Here, $\Gamma_{\text{th}} = 5$ dB, and the null spaces represent regions where SINR $< \Gamma_{\text{th}}$.}
	
	\label{cf}
\end{figure*}
\begin{table}[t]
	\centering
	\caption{Simulation parameters.}
	\label{tab:setup}
	\begin{tabular}{lc}
		\toprule
		\textbf{Parameter} & \textbf{Value} \\
		\midrule
		Room $(L{\times}W{\times}H)$ & $8{\times}8{\times}3$ m \\
		Grid resolution $\Delta$ & $0.2$ m \\
		Transmitters $N_{\text{TX}}$ & $225$ ($15{\times}15$) \\
		Power per VCSEL $P_t$ & $10$ mW \\
		Detector area $A_d$ & $10^{-4}$ m$^2$ \\
		FoV / index $(\Psi_{\text{FoV}}, n)$ & $75^\circ$, $1.5$ \\
		RL $(\alpha,\gamma,N_{\text{ep}})$ & $0.1$, $0.9$, $2000$ \\
		RL $(\epsilon_{\max}, \epsilon_{\min}, \lambda)$ & $1$, $0.01$, $0.005$ \\
		\bottomrule
	\end{tabular}
	\label{parameters}
\end{table}

We evaluate the proposed approach in a simulated indoor environment using the model described in Section~\ref{sys-model}. The simulation parameters are listed in Table \ref{parameters}. \textcolor{black}{We also note that, the analysis in this paper focuses on a controlled single-access-point indoor scenario with a fixed receiver plane and quasi-static user height states. This setup is intentionally chosen to isolate the impact of beam divergence control and learning-based decision making on coverage performance. Dynamic effects such as user mobility, blockage, and multi-user interference are not considered in the present study and are left for future work.}

First, to obtain deeper insight into the importance of optimizing the beam divergence angle, Fig.~\ref{cf} illustrates the simulated three-dimensional coverage performance for a fixed VCSEL divergence angle of $\theta_{\text{divergence}} = 15^{\circ}$. The figure clearly exposes the fundamental limitation of static beam configurations in indoor OWC environments. As the receiver height increases from 0.5~m to 2.0~m, the overall coverage quality deteriorates sharply from approximately 80\% down to 48\%. This degradation arises from of using a fixed divergence angle. At lower heights ($h_r = 0.5$~m), the $15^{\circ}$ beam is excessively wide, producing strong inter-beam interference and scattered coverage holes. Conversely, at higher heights ($h_r = 2.0$~m), the same divergence  cannot provide adequate spatial illumination, resulting in power loss and incomplete coverage across the receiver plane. These results demonstrate that the trade-off between beam power concentration and spatial spread is highly sensitive to the vertical geometry of the link. The optimal divergence angle changes dynamically with receiver height, but this dependence is nonlinear, non-convex, and analytically intractable due to the discrete nature of the transmitter array and the non-differentiable coverage function. 
%
To address this limitation, our proposed RL–based control framework enable autonomous and adaptive adjustment of the divergence angle.

\begin{figure}
	\begin{center} 
			\includegraphics[width=3.2 in]{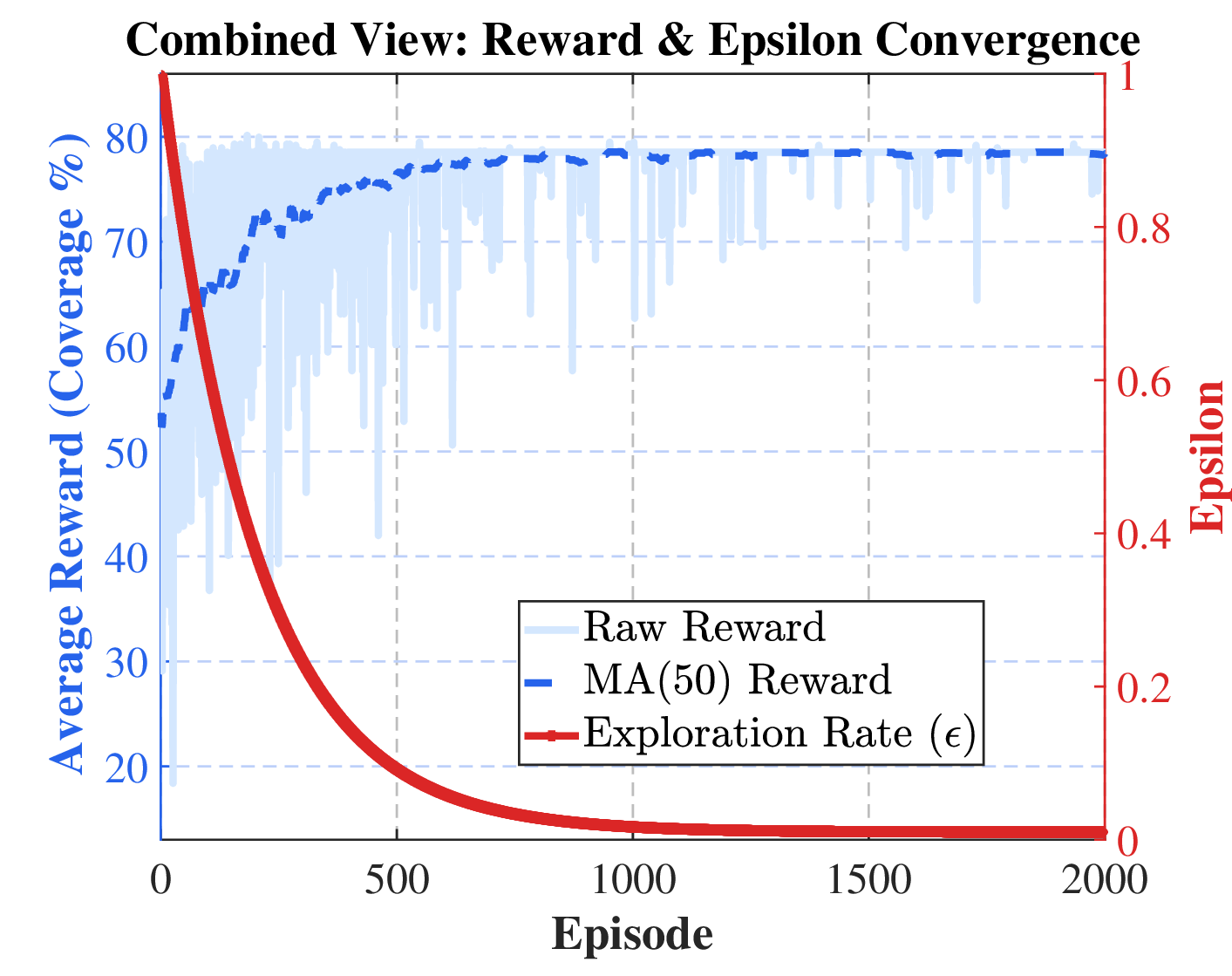}
			\caption{Convergence behavior of the proposed RL-based divergence control framework. The blue curve shows the evolution of the average coverage reward per episode (with a 50-episode moving average (MA)), while the red curve represents the exploration rate ($\epsilon$).}
			\label{rl_convergence}
		\end{center}
\end{figure}
%
\subsection{Convergence Behavior and Learning Stability}
The convergence behavior of the proposed RL-based divergence control framework is shown in Fig.~\ref{rl_convergence}. The figure jointly depicts the average reward, representing the achieved coverage percentage, and the exploration rate $\epsilon$ as functions of the training episode. The agent begins with a fully exploratory policy ($\epsilon = 1$), leading to large fluctuations in the early-stage reward due to random exploration of suboptimal actions. As training progresses, the exploration rate decays exponentially, thereby allowing the agent to gradually exploit the accumulated experience and reinforce high-reward actions.

Correspondingly, the moving average of the coverage reward exhibits a clear upward trend, converging smoothly to a stable plateau around 85\%. This indicates that the agent has successfully learned a consistent mapping between the environmental state (receiver height) and the optimal divergence angle $\theta_{\text{divergence}}$. The absence of significant oscillations in the steady-state phase further demonstrates the stability and robustness of the learned policy. It should be noted, however, that the achieved coverage does not reach 100\% in the considered setup. This limitation arises not from the learning framework itself but from the physical constraints of the considered system model, where each user location is illuminated by a single VCSEL beam. To further enhance coverage, particularly in areas with overlapping or weak illumination, multiple-beam support (i.e., assigning more than one transmitter per user) would be required \cite{safi3d}. Such configurations can provide redundancy and spatial diversity, thereby mitigating coverage gaps and pushing system performance closer to full spatial utilization.
\vspace{-3mm}
\subsection{Results and Learned Divergence Policy }

The convergence profile underscores two key observations. First, the RL agent effectively balances exploration and exploitation through the adaptive decay of $\epsilon$, facilitating rapid policy refinement despite the environment’s non-differentiable nature. Second, the resulting policy attains near-optimal performance comparable to exhaustive optimization methods, yet with substantially lower computational complexity and without relying on prior system knowledge. These findings validate the efficacy of the proposed RL framework in autonomously optimizing beam divergence for dynamic three-dimensional OWC systems.

To further confirm the effectiveness of the learned policy and assess its physical consistency, Table~\ref{tab:rl_policy} summarizes the optimal divergence angles and corresponding coverage percentages obtained across different receiver heights. From the table, it can be realized that the RL agent successfully captures the nonlinear dependence between the optimal divergence angle and the receiver height to adaptively balance beam spread and interference. Accordingly, at lower heights ($h_r = 0.5$–1.0~m), the policy selects relatively narrow beam angles ($8^{\circ}$–$9^{\circ}$), effectively concentrating optical power and minimizing inter-beam overlap. As the receiver height increases, the policy learns to gradually widen the divergence angle, reaching $14^{\circ}$ at $h_r = 2.0$~m, which offers geometric beam expansion and maintains sufficient illumination across the receiver plane.

Quantitatively, the learned policy achieves coverage levels closely matching those obtained through exhaustive optimization, with deviations of less than 3\% at most heights. For example, at $h_r = 0.5$~m and $h_r = 1.0$~m, the learned divergences of $9^{\circ}$ and $8^{\circ}$ yield coverage of 92.39\% and 85.72\%, respectively, nearly identical to the exhaustive-search baselines of 93.33\% and 88.81\%. Even at the most challenging configuration ($h_r = 2.0$~m), the RL policy maintains 50.03\% coverage compared to 53.8\% from the exhaustive search (ES), reflecting the inherent geometric limitations of single-beam operation rather than a limitation of the learning framework. \textcolor{black}{The observed saturation in coverage probability at higher receiver heights (approximately 50\% at 2 m) is primarily a consequence of the single-beam-per-user transmission model rather than a limitation of the learning framework itself. These results motivate future extensions toward multi-beam or cooperative transmission strategies, where multiple VCSEL elements jointly serve a user or region.}

\textcolor{black}{Although the proposed Q-learning formulation operates on discretised receiver heights and a finite set of beam divergence actions limits direct generalisation to fully continuous or rapidly time-varying environments, it enables stable learning, low computational overhead, and clear physical interpretability. Extensions to continuous state spaces using function approximation or deep reinforcement learning are a natural next step, but were intentionally avoided here to preserve analytical transparency and avoid overfitting in small-sample regimes.}

Beyond accuracy, the proposed RL approach offers a  computational advantage critical for scaling to real-world environments and mitigating {environmental drift}. While the online inference complexity for both the ES and the RL policy is technically $\mathcal{O}(1)$ (a single table lookup), this equivalence neglects the impractical {offline cost} required to initially generate the ES table, which scales prohibitively as $\mathcal{O}(|\mathcal{A}| N_{\text{grid}} N_{\text{TX}})$. This high cost must be paid every time the environment changes (e.g., if a new floor plan or large furniture alters the effective optical channel). In contrast, the RL framework shifts the immense computational burden into a single, efficient, data-driven training process. The RL training complexity scales only with the number of episodes and states ($\mathcal{O}(N_e  |\mathcal{S}|)$) and is {independent of the expensive spatial sampling} ($N_{\text{grid}}$ and $N_{\text{TX}}$). This minimal memory and training dependency on the physical space prove that the RL framework is the only scalable and maintainable solution for dense VCSEL arrays and high-speed operation.

\begin{table}[t]
	\centering
	\small
	\caption{Learned Optimal Divergence Angles and Coverage Performance Across Receiver Heights}
	\label{tab:rl_policy}
	\renewcommand{\arraystretch}{1.1}
	\begin{tabular}{c|cc|cc}
		\hline
		\textbf{$h_r$ [m]} & \multicolumn{2}{c|}{\textbf{RL Learned Policy}} & \multicolumn{2}{c}{\textbf{Exhaustive Baseline}} \\
		\cline{2-5}
		& $\theta_{\text{divergence}}$ [$^\circ$] & Cov. [\%] & $\theta_{\text{divergence}}$ [$^\circ$] & Cov. [\%] \\
		\hline
		0.5 & 9  & 92.39 & 10 & 93.33 \\
		1.0 & 8  & 85.72 & 10 & 88.81 \\
		1.5 & 9  & 85.96 & 8  & 89.70 \\
		2.0 & 14 & 50.03 & 11 & 53.80 \\
		\hline
	\end{tabular}
\end{table}

\begin{table}[t]
	\centering
	\caption{Computational Complexity Analysis of ES and RL-Based Divergence Control}
	\label{tab:complexity}
	\renewcommand{\arraystretch}{1.1}
	\small
	\begin{tabular}{l|c|c}
		\hline
		\textbf{Method} & \textbf{Offline} & \textbf{Online} \\
		\hline
		\textbf{ES} & $\mathcal{O}(|\mathcal{A}| N_{\text{grid}} N_{\text{TX}})$ & $\mathcal{O}(1)$ \\
		\textbf{RL (Q-Learn)} & $\mathcal{O}(N_e \cdot |\mathcal{S}|)$ & $\mathcal{O}(1)$ \\
		\hline
	\end{tabular}
	
	\vspace{1mm}
	\raggedright
	\footnotesize
	Both methods require $\mathcal{O}(|\mathcal{S}| \times |\mathcal{A}|)$ memory. $N_e$: training episodes; $N_{\text{grid}}$: grid points. Large $N_{\text{grid}}$ and $N_{\text{TX}}$ make ES offline cost prohibitive.
\end{table}
\section{Conclusion}
\label{Conclusion}
WE presented an RL framework for adaptive beam divergence control in 3D indoor OWC systems. The RL agent learns an efficient mapping between receiver height and beam divergence without relying on analytical models or ES. The resulting lightweight Q-table policy offers near-optimal coverage and real-time adaptability. Results showed that data-driven divergence control enhances scalability and responsiveness in dense VCSEL-based networks. \textcolor{black}{Future work will consider dynamic user mobility, blockage effects, and multi-user interference, as well as continuous-state learning approaches and cooperative multi-beam transmission. These extensions are expected to further improve robustness and scalability in realistic indoor deployments.}
\section*{Acknowledgement}
{This work is a contribution by Project REASON, a UK Government funded project under the Future Open Networks Research Challenge (FONRC) sponsored by the Department of Science Innovation and Technology (DSIT). The authors also acknowledge support by the Engineering and Physical Sciences Research Council (EPSRC) under grants EP/X04047X/1 - EP/Y037243/1 ``Platform for Driving Ultimate Connectivity (TITAN).''}
	\bibliographystyle{IEEEtraN}


\end{document}